\documentclass[a4paper, 11pt, onecolumn, oneside]{article}
\usepackage[T1]{fontenc}
\usepackage[utf8]{inputenc}
\usepackage{makeidx}
\usepackage{graphicx}
\usepackage{float}
\usepackage{amssymb}
\usepackage{amsmath}
\usepackage{latexsym}
\usepackage{wasysym}
\usepackage{amsthm}
\usepackage{natbib}
\usepackage{}
\usepackage{amsfonts}
\usepackage{amsmath}
\usepackage{amssymb}
\usepackage{amstext}

\begin{document}

\title{Generalized Extreme Value distribution parameters as dynamical indicators of Stability}

\author{Faranda, Davide\\
\small{\textit{Department of Mathematics and Statistics, University of Reading;}}\\
\small{\textit{Whiteknights, PO Box 220, Reading RG6 6AX, UK.} d.faranda@pgr.reading.ac.uk}\\ \\
Lucarini, Valerio\\
\small{\textit{Klimacampus, Universit\"at Hamburg;}}\\
\small{\textit{Grindelberg 5, 20144, Hamburg, Germany;}}\\
\small{\textit{Department of Mathematics and Statistics, University of Reading;}}\\
\small{\textit{Whiteknights, PO Box 220, Reading RG6 6AX, UK.} v.lucarini@reading.ac.uk}\\ \\
Turchetti, Giorgio\\
\small{\textit{Department of Physics, University of Bologna.INFN-Bologna}}\\
\small{\textit{Via Irnerio 46, Bologna, 40126, Italy.} turchett@bo.infn.it}\\ \\
Vaienti, Sandro\\
\small{\textit{UMR-6207, Centre de Physique Th\'eorique, CNRS, Universit\'es d'Aix-Marseille I,II,}}\\
\small{\textit{Universit\'e du Sud Toulon-Var and FRUMAM}}\\
\small{\textit{(F\'ed\'eration de Recherche des Unit\'es de Math\'ematiques de Marseille);}}\\
\small{\textit{CPT, Luminy, Case 907, 13288 Marseille Cedex 09, France.}}\\
\small{vaienti@cpt.univ-mrs.fr}\\
}
\date{}
\maketitle

\begin{abstract}
We introduce a new dynamical indicator of stability based on the Extreme Value statistics showing that it provides an insight on the local stability properties of dynamical systems. The indicator perform faster than other based on the iteration of the tangent map since it requires only
the evolution of the original systems and, in the chaotic regions, gives further information about the information dimension of the attractor. A numerical validation of the method is presented through the analysis of the motions in a Standard map.

\end{abstract}

\section{Introduction}
The analysis of stability for discrete and continuous time dynamical systems is of fundamental importance to get insights in the dynamical structure of a system. The distinction among regular and chaotic orbits can be easily made in a dissipative case whereas for conservative systems it is usually an hard task especially when we are dealing with many degrees of freedom or when the phase space is splitted into chaotic regions and regular islands. A large number of tools known as indicators of stability have been developed to accomplish this task: Lyapunov Characteristic Exponents (LCEs)    
\citep{wolf1985determining}, \citep{rosenstein1993practical}, \citep{skokos2010lyapunov}
and the indicators related to the Return Time Statistics \citep{kac1934notion}, \citep{gao1999recurrence}, \citep{hu2004statistics} , \citep{buric2005statistics} have been used from a long time as efficient indicators. Nevertheless, in the recent past, the need for computing stability properties with faster algorithms and for systems  with many degrees of freedom resulted in a renewed interest in the technique and different dynamical indicators have been introduced. The Smaller Alignment Index  (SALI) described in \citet{skokos2002smaller} and \citet{skokos2004detecting} ,the Generalized Alignment Index (GALI), introduced in \cite{skokos2007geometrical}and the Mean Exponential Growth factor of 
Nearby Orbits (MEGNO) discussed in \citet{cincotta2003phase},\citet{gozdziewski2001global} are suitable to analyse the properties of a single orbit. They are based on the divergence of nearby trajectories and require in principle the knowledge of the exact dynamics. Another class of indicators is based on the round off error properties and has been discussed in \citet{faranda2011reversibility}: the divergence between two trajectories starting from the same initial condition but computed with different numeric precision can be used to overlight the dynamical structure. The so called Reversibility Error that measure the distance between a certain initial condition and the end point of a trajectory iterated forward and backward for the same number of time steps give basically the same informations.\\
Other tools such as Fidelity and Correlation decay can be  successfully used to 
characterize stability properties of ensemble of orbits as explained in \citet{vaienti2007random} and \citet{turchetti2010relaxation}. Eventually,  
Frequency Map Analysis has been used to provide informations about the resonance structure of a system \citep{laskar1999introduction}, \citep{robutel2001frequency}.\\

Although these indicators have been developed to accomplish the same task, each of them presents specific features and limitations and often it is necessary to combine a certain number of indicators to get quantitative information about the dynamics. Another aspect to consider when comparing the efficiency of chaos indicators is the computational cost: each
variational method  (SALI, MEGNO, mLCE) needs to iterate both the map and the tangent map forward during $n$ steps. The latter
map, although linear, is the computationally most expensive due to the fact that it needs the evaluation
of the Jacobian matrix at every step. On the other hand the Round off error methods are the computationally less expensive as they require only the iteration of the dynamics. Fidelity and Correlations are usually computed using Monte Carlo simulations and therefore they are inaccessible for systems with many degrees of freedom \cite{turchetti2010relaxation}.
Therefore, it is clear that there is still the need to introduce versatile indicators that distinguish regular from chaotic behaviors and possibly give further informations on the dynamics. The pourpose of this paper is to use the series of extrema of some specific observable computed using the orbits of dynamical systems as a new indicator of stability.\\

Extreme Value Theory was originally introduced by \citet{fisher}, \citet{gnedenko}to study the maxima of a series of independent and identical distributed variables: under very general hypothesis a limiting distribution called Generalised Extreme Value (GEV) distribution exists for the series of extremes. An extensive account of recent results and relevant applications is given in \citet{ghil2010extreme}. In the recent past this theory has been adapted to study the output of dynamical systems.  As we will explain in detail in the next section it is not trivial to observe the asymptotic GEV distribution in dynamical systems: some sort of independence of maxima must be recovered by requiring certain mixing conditions on the orbits. Furthermore, we need to introduce some peculiar observables that satisfy the condition proposed by Gnedenko on the parent distribution of data: they are related to the closest return of a trajectory in a ball centered around the starting point and therefore allow a very detailed description of the dynamics in the neighborhood of the initial condition.\\
If all these requirement are satisfied it is possible to observe an Extreme Value statistics that converges to the GEV distribution family. The parameters of the distribution are intimately related with some relevant dynamical quantities such as the dimension of the attractor \citep{freitas}, \citep{faranda2011extreme}. We will use these features to show the reliability of GEV distribution parameters in discriminate regular from chaotic behaviors pointing out the further informations that is possible to extract regarding the dynamics.\\
The paper is organised as follows: in section 2 we explain how to introduce an EV statistics in dynamical systems pointing out the difference between regular and chaotic orbits. In section 3 we describe the numerical algorithm and procedure used to compute the parameters of the GEV distribution. Eventually, in section 4 we perform some tests on the Standard Map to validate numerically the use of GEV parameters as dynamical indicators.

\section{Extreme Statistics as dynamical indicator}

\subsection{Extreme Value Theory in dynamical systems}

\citet{gnedenko}  studied the convergence of maxima of i.i.d.
variables $$X_0, X_1, .. X_{m-1}$$ with cumulative distribution
(cdf) $F(x)$ of the form:
$$F(x)=P\{a_m(M_m-b_m) \leq x\}$$
Where $a_m$ and $b_m$ are normalizing sequences and $M_m=\max\{ X_0,X_1, ..., X_{m-1}\}$. It may be rewritten as $F(u_m)=P\{M_m \leq u_m\}$ where $u_m=x/a_m +b_m$. Under general hypothesis on the nature of the parent distribution of data, \citet{gnedenko} show that the distribution of maxima, up to an affine change of variable,  obeys to one of the following three laws:
\begin{itemize}
\item Type 1 ({\em Gumbel}). \begin{equation}E(x)=\exp(-e^{-x}), \ -\infty<x<\infty \label{gum}\end{equation}
\item Type 2 ({\em Fr\'echet}). \begin{equation}E(x)= \begin{cases}
 0,\  x\le 0 \\
\exp(-x^{-\xi}),\  \mbox{for some} \ \xi>0, \ x>0
\end{cases} \label{fr} \end{equation}

\item Type 3 ({\em Weibull}). \begin{equation} E(x)= \begin{cases}
\exp(-(-x)^{-\xi}),\  \mbox{for some} \  \xi>0, \ x\le 0 \\
1,\  x>0
\end{cases} \label{wei}\end{equation}

\end{itemize}

Let us define the right endpoint $x_F$ of a distribution function $F(x)$ as:
\begin{equation}
x_F=\sup\{ x: F(x)<1\}
\end{equation}

then, it is possible to compute normalizing sequences   $a_m$ and $b_m$ using the following corollary of Gnedenko's theorem :\\
\textbf{Corollary (Gnedenko):}  \textit{The normalizing sequences $a_m$ and $b_m$ in the convergence of normalized maxima $P\{a_m(M_m - b_m) \leq x\} \to F(x)$ may be taken (in order of increasing complexity) as:}

\begin{itemize}

\item \textit{Type 1:} $\quad a_m=[G(\gamma_m)]^{-1}, \quad b_m=\gamma_m$;

\item \textit{Type 2:} $\quad a_m=\gamma_m^{-1}, \quad b_m=0$;

\item \textit{Type 3:} $\quad a_m=(x_F-\gamma_m)^{-1}, \quad b_m=x_F$;

\end{itemize}
\textit{where}
\begin{equation}
\gamma_m=F^{-1}(1-1/m)=\inf\{x; F(x) \geq 1-1/m\}
\label{gamma}
\end{equation}

\begin{equation}
G(t)=\int_t^{x_F} \frac{1-F(u)}{1-F(t)}du, \quad  t<x_F
\label{gneg}
\end{equation}

In \citep{faranda2011numerical} we have shown that this approach is equivalent  to fit unnormalized data directly  to a single family of generalized distribution called GEV distribution with cdf:

\begin{equation}
F_{G}(x; \mu, \sigma,
\xi')=\exp\left\{-\left[1+{\xi'}\left(\frac{x-\mu}{\sigma}\right)\right]^{-1/{\xi'}}\right\}
\label{cumul}
\end{equation}

which holds for $1+{\xi'}(x-\mu)/\sigma>0 $, using $\mu \in
\mathbb{R}$ (location parameter) and $\sigma>0$ (scale parameter) as
scaling constants in place of $b_m$, and $a_m$ \citep{pickands}, in
particular, in \citet{faranda2011numerical} we have shown that the
following relations hold:

$$\mu=b_m \qquad \sigma=\frac{1}{a_m}. $$

${\xi'} \in \mathbb{R}$ is the shape parameter also called the tail index: when ${\xi'} \to 0$, the distribution corresponds to a
Gumbel type ( Type 1 distribution).  When the index is positive, it corresponds to a Fr\'echet (Type 2 distribution); when the index is negative, it corresponds to a Weibull (Type 3 distribution).\\

In the last decade many works focused on the possibility of treating time series of observables of deterministic dynamical system  using EVT.  For example,  \citet{balakrishnan} and more recently \citet{nicolis} and \citet{haiman} have shown that for regular orbits of dynamical systems we don't expect to find convergence to EV distribution.\\
The first rigorous  mathematical approach to extreme value theory in dynamical systems goes back to the pioneer paper by P. Collet in 2001 \citep{collet2001statistics}. Important contributions have successively been given by  \citet{freitas2008}, \citet{freitas},\citet{freitas2010extremal} and by \citet{gupta2009extreme}. The starting point of all these investigations was
to associate to the stationary stochastic process given by the dynamical
 system, a new stationary independent sequence which enjoyed one of the
  classical three extreme value laws, and this laws could be pulled back
  to the original dynamical sequence. 
  
  Let us consider a dynamical systems $(\Omega, {\cal B}, \nu, f)$, where $\Omega$ is the
    invariant set in some manifold, usually $\mathbb{R}^d$, ${\cal B}$ is the
     Borel $\sigma$-algebra, $f:\Omega\rightarrow \Omega$ is a measurable map
      and $\nu$ a  probability $f$-invariant Borel measure.\\
In order to adapt the extreme value theory to dynamical systems, we will consider the stationary stochastic process $X_0,X_1,...$  given by:

\begin{equation}
X_m(x)=g(\mbox{dist}(f^m (x), \zeta)) \qquad \forall m \in \mathbb{N}
\label{sss}
\end{equation}

where 'dist' is a distance on the ambient space  $\Omega$, $\zeta$ is a given point and $g$ is an observable function, and whose partial maximum is defined as:

\begin{equation}
{M_m}= \max\{ X_0, ... , X_{m-1} \}
\label{maxi}
\end{equation}

The probability measure will be here the invariant measure $\nu$ for the dynamical system. We will also suppose that our systems which verify the condition $D_2$ and $D'$ which will allow us to use the EVT for i.i.d. sequences. As we said above, we will use   three types of observables $g_i,i=1,2,3$,  suitable to obtain one of the three types of EV distribution  for normalized maxima:

\begin{equation}
g_1(x)= -\log(\mbox{dist}(x,\zeta))
\label{g1}
\end{equation}

\begin{equation}
g_2(x)=\mbox{dist}(x, \zeta)^{-1/\alpha}
\label{g2}
\end{equation}

\begin{equation}
g_3(x)=C - \mbox{dist}(x,\zeta)^{1/\alpha}
\label{g3}
\end{equation}

where $C$ is a constant and $\alpha>0 \in \mathbb{R}$.\\

Using these observables we can obtain  convergence  to the Type 1,2,3 distribution if  one can prove  two sufficient conditions called $D_2$ and $D'$ and which we briefly explain here: these conditions  basically require a sort of independence of the stochastic dynamical sequence in terms of uniform mixing condition on the distribution functions. In particular condition $D_2$, introduced in its actual  form by Freitas-Freitas \citet{freitas2008},   could be checked directly by estimating the rate of decay of correlations for H\"older observables.\\
If $X_m, m\ge 0$ is our stochastic process, we can define  $M_{j,l}\equiv \{X_j, X_{j+1},\cdots,X_{j+l}\}$ and we put $M_{0,m}=M_m$.\\
The condition $D_2(u_m)$ holds for the sequence $X_m$ if for any integer $l, t,m$ we have $$|\nu (X_0>u_m, M_{t,l}\le u_m)-\nu (X_0>u_m)\nu (M_{t,l}\le u_m)|\le \gamma(m,t)$$ where $\gamma(m,t)$ is non-increasing in $t$ for each $m$ and $m\gamma(m,t_m)\rightarrow 0$ as $m\rightarrow \infty$ for some sequence $t_m=o(m)$, $t_m \rightarrow \infty$.\\ 
We say condition $D'( u_m)$ holds for the sequence $X_m$ if $$\lim_{l\rightarrow \infty}\limsup_{m} m\sum_{j=1}^{[m/l]}\nu(X_0>u_m, X_j>u_m)=0$$ \\
Instead of checking the previous conditions, we can use another results that established a connection between the extreme value laws and the statistics of first return and hitting times, see the papers by \citet{freitas} and \citet{freitasNuovo}. They showed in particular that for dynamical systems preserving an absolutely continuous invariant measure or a singular continuous invariant measure $\nu$, the existence of an exponential hitting time statistics on balls around $\nu$ almost any point $\zeta$ implies the existence of extreme value laws for one of the observables of type $g_i, i=1,2,3$ described above. The converse is also true, namely if we have an extreme value law which applies to the observables of type $g_i, i=1,2,3$ achieving a maximum at $\zeta$, then we have exponential hitting time statistics to balls with center $\zeta$. Recently these results have been generalized to local returns around balls centered at periodic points \citep{freitas2010extremal}.\\

\subsection{Extreme Value Statistics in mixing and regular maps}

In \citet{faranda2011numerical} and  \citet{faranda2011extreme} we have analised both from an analytical and numerical point of view the Extreme Value distribution in a wide class of low dimensional maps showing that, when the conditions $D'$ and $D_2$ are verified, the block maxima approach can be used to study extrema. It consists of dividing the data series of length $k$ of some observable into $n$ bins each containing the same number $m$ of observations, and selecting the maximum (or the minimum) value in each of them \citep{coles}.
Using $g_i$ observable functions we have shown that a first order approximation of the GEV parameters in mixing maps can be written in terms of $m$ (or equivalently $n$) and the dimension of the attractor $d$:

For $g_1$ type observable:

\begin{equation}
\sigma= \frac{1}{d} \qquad \mu \sim \frac{1}{d}\ln(k/n) \qquad \xi'=0
\label{g1res}
\end{equation}

For $g_2$ type observable:

\begin{equation}
\sigma\sim n^{-1/(\alpha d)} \qquad \mu \sim n^{-1/(\alpha d)} \qquad \xi'=\frac{1}{\alpha d}
\label{g2res}
\end{equation}

For $g_3$ type observable:

\begin{equation}
\sigma\sim n^{1/(\alpha d)} \qquad \mu = C \qquad \xi'=\frac{1}{\alpha d}
\label{g3res}
\end{equation}

while the higher order terms contain explicit dependence on the density measure.\\

For regular maps independently
on the observable chosen, for periodic or quasi-periodic orbit we do not observe convergence to the GEV distribution.
In \citet{faranda2011numerical} we have analysed what happens when observable $g_i$ are used to study empirical maxima distribution 
in regular maps: \citet{nicolis}  have shown how it is possible to obtain an analytical EV distribution which does not belong to GEV family choosing a simple observable: they considered the series  of  distances  between the iterated trajectory and  the initial condition:

$$ Y_m(x=f^t\zeta)= \mbox{dist}(f^t\zeta, \zeta) \qquad \hat{M}_m=\min \{Y_0, ... Y_{m-1} \}$$

For this observable they have shown that the cumulative distribution $F(x)=P\{a_m(\hat{M}_m-b_m) \leq x\}$  of a uniform quasi periodic motions is not smooth but piecewise linear (\citet{nicolis}, Figure 3).
$F(x)$ must correspond  to a density distribution  continuous obtained as a composition of Heaviside step functions:  each box must be related to a change in the slope of $F(x)$. \\
By applying the observable $g_i$ we just remodulate this piecewise linear $F(x)$ but it is clear that we don't obtain any kind of convergence to the GEV distribution. In terms of density functions, we can observe multimodal distributions and the number of modes and their positions are highly dependent on both $n,m$ and initial conditions \citet{faranda2011numerical}. \\
Using \citet{nicolis} results it is possible to understand the empirical distribution obtained once applied $g_i$ functions: since density distribution  of $\hat{M}_m$ is  a composition of box functions, when we apply $g_i$ observables we  modulate it changing the shape of the boxes. Therefore, we obtain a multi modal distribution modified according to the observable functions $g_i$. In the case of pure periodic motion this kind of extreme value distribution must asymptotically be a Dirac delta as we pick up always the same $Y_i$.\\

\section{The numerical algorithm}

As we have already said  to observe a GEV distribution of maxima  orbits must satisfie $D'$ and $D_2$ conditions or an exponential decay of the Hitting Time Statistics whereas for periodic or quasi-periodic motions we have different distributions but never a GEV. This give us a way to discriminate the kind of motion simply looking at the extreme value statistics. From a practical point of view we can introduce a simple algorithm to perform this task:

\begin{enumerate}
\item Compute the orbit of the dynamical system for $k$ iterations.
\item Compute the series $X_m(x)=g(\mbox{dist}(f^m (x), \zeta))$ where $\zeta$ is the initial condition.
\item Divide the series in $n$ bins each containing $m$ data.
\item Take the maximum in each bin and test if this empirical distribution give the parameters expected by the theory.
\end{enumerate}

In the next section we will describe how to use it operationally in a meaningful example of a dynamical systems that present
coexistence of regular and chaotic motions and therefore allow to test the validity of our indicator: the Standard map.
Before presenting the results, we need to clarify the numerical inference
procedure that we use to obtain the parameters of the GEV distribution.
In \citet{faranda2011numerical}  we have
used a Maximum Likelihood Estimation   (MLE) procedure working both on pdf and cdf (cumulative distribution
function), since our distributions were absolutely continuous and
the minimization procedure was well defined. In a general case, when we are dealing also with singular measures,  we may have a cdf which is not anymore absolutely continuous  and consequently the fitting procedure via MLE
could give wrong results \citep{faranda2011extreme}. To avoid these problems
  we have used an L-moments estimation
   as detailed in \citet{hosking1990moments}.
   This procedure is completely discrete and can be used both for
   absolutely  continuous or  singular continuous  cdf. The L-moments
   are summary statistics for probability distributions and data samples.
    They are analogous to ordinary moments which meant that
    they provide measures of location, dispersion, skewness,
     kurtosis, but are computed from linear combinations of
     the data values, arranged in increasing order (hence the prefix L). Asymptotic approximations to sampling distributions are better for L-moments than for ordinary moments [\citet{hosking1990moments}, Figure 4].  The relationship between the moments and the parameters of the GEV distribution are described in \citet{hosking1990moments}, while the 95\% confidence intervals has been derived using a bootstrap procedure.\\
We have explained in \citet{faranda2011extreme} the issues of applying directly a Kolmogorov Smirnov or a Chi-square test to see if the data really belongs to the GEV distribution when the measure is singular: even if both tests fail the GEV model is reasonable as it is the closest continuous representation of the empirical discrete distribution we get.\\
Therefore, instead of using as a dynamical indicator the goodness of the fit to the continuous GEV model, we will check the deviations of the parameters with respect to the theoretical expected values.\\

Our results will be studied against the Divergence of two nearby trajectories and the Reversibility error. In \citet{faranda2011reversibility} we have shown that these indicators give insights into the structure of a dynamical system as well as others such as SALI, MEGNO and mLCE with which these indicators have been compared.\\
We briefly recall here the definitions and we remind to \citet{faranda2011reversibility} for further clarifications.

The arithmetic operations such as 
sums or multiplications 
imply a round  off, which  propagates  the error affecting each number. Round off algebraic
procedures are hardware dependent as detailed in \citet{knuth1973art}.
Unlike the case of stochastic perturbations, the error strongly depends on $x$.
Suppose we are given a map $f^t(x)$ then we will indicate with $f^t_*(x)$ the correspondent numerical map both evaluated at the $t$-th iteration.
The Divergence of orbits is defined as:
\begin{equation}
  \Delta_t=\mbox{dist}(f_S^t(x),f_D^t(x)),
   \label{eq:def_diver}
\end{equation}
where $f_S^t$ and $f_D^t$ stand for single and double precision iterations respectively, and $\mbox{dist}$ is a suitable metrics.

If the map is invertible we can also define the Reversibility error as

\begin{equation}
  R_t=\mbox{dist}(f_*^{-t}\circ f_*^t (x), x)
  \label{eq:def_rever}
\end{equation}

which is non zero since the numerical inverse $f_*^{-1}$  of the map 
is not exactly the inverse of $f_*$ namely $f_*^{-1}\circ f_*(x)\not = x$.
Obviously the reversibility error is much easier to compute than the divergence 
of orbits (if we know explicitly the inverse map) and the information it provides 
is basically the same as the latter. 
Both quantities give an average linear growth for a regular map together with an 
exponential growth for a chaotic map having positive Lyapounov exponents and strong mixing properties. When computing $R_t$  we will set $f_* = f_S$
in order to compare with $\Delta_t$.\\

\section{A case of study: the Standard Map}

The Standard map (also known as Chirikov-Taylor map or Chirikov standard map) is an area-preserving chaotic map defined on the bidimensional torus. It can be thought as   a stick that is free of the gravitational force, which can rotate frictionless in a plane around an axis located in one of its tips, and which is periodically kicked on the other tip. This mechanical system  is usually called a kicked rotator. It is defined by:

\begin{equation}
\begin{cases}
y_{t+1}= y_t - \frac{K}{2\pi}\sin(2 \pi x_t) & \mod 1  \\
x_{t+1}= x_t + y_t+1 & \mod 1  \\ 
\end{cases}
\label{stdmap}
\end{equation}

\begin{figure}[ht!]
  \centering  
  \includegraphics[width=0.95\textwidth]{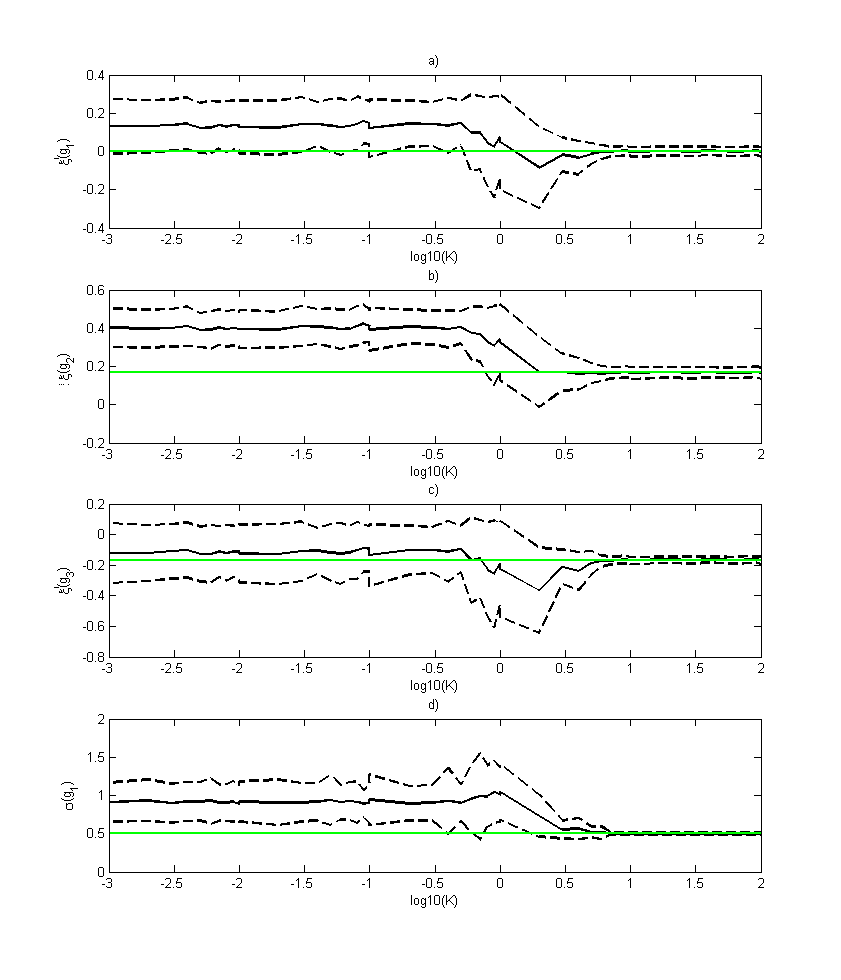}
  \caption{Standard map: GEV parameters averaged over 500 different initial condition centered in $x_0=0.305, y_0=0.7340$ VS $K$. a) $\xi'(g_1)$, b) $\xi'(g_2)$, c) $\xi'(g_3)$, d) $\sigma(g_1)$. }
  \label{perk}
\end{figure}

Standard map is one of the most widely-studied examples of dynamical chaos in physics.
It can be regular or chaotic, depending on the strength of the impulses: stronger kicks lead to chaotic behaviors.  The variables $y$ and $x$ respectively represent the angular position and angular momentum of the stick at the $t$-th kick.\\
For $K<<1$ the motion follows quasi periodic orbits for all initial conditions, whereas if $K>>1$ the motion turns to be chaotic and irregular.
An interesting behavior is achieved when $K \sim 1$: in this case we have coexistence of regular and chaotic motions depending on the initial
conditions chosen.\\
First of all we have taken an ensemble of 500 initial conditions centered around $x_0, y_0$ in a small subset of the bidimensional torus. After iterating the map for $k=10^6$ iterations, we have selected $n=1000$ maxima each in $m=1000$ observations. These values are in agreement with what emerged from the numerical study presented in \citet{faranda2011numerical}. Once computed the parameters of GEV distribution for each realization, we have averaged them over the different initial conditions  checking their convergence towards the expected theoretical values when $K$ is varied from $K=10^{-4}$ up to $K=10^2$.
As indicators, looking at the equations \ref{g1res}-\ref{g3res}, we have selected the shape parameters for the three type observables $\xi'(g_1), \xi'(g_2), \xi'(g_3)$ and the scale parameter for the type 1 observable $\sigma(g_1)$. From all the parameters set, these are the ones that have a dependence on the dimension of the attractor $d$ but do not have a dependence of $m$, therefore, once we are in the asymptotic regime, the results are independent on the number of observations in each bin. Results are presented in figure \ref{perk}. In this example we set  $\alpha=3$ for  $g_2$ and $g_3$ observable. For each parameters, the averaged value is represented with a solid line whereas the dotted lines represent one standard deviations of the ensemble. It is clear that for $K>1$ the parameters converge towards the theoretical values whereas for regular motions, the computed parameters are not representative of a GEV distribution and exhibit a spread that is more than five times bigger with respect to the chaotic counterpart. The results are similar if we change  the initial conditions and the value of $\alpha$.\\

\begin{figure}[ht!]
  \includegraphics[width=1.0\textwidth]{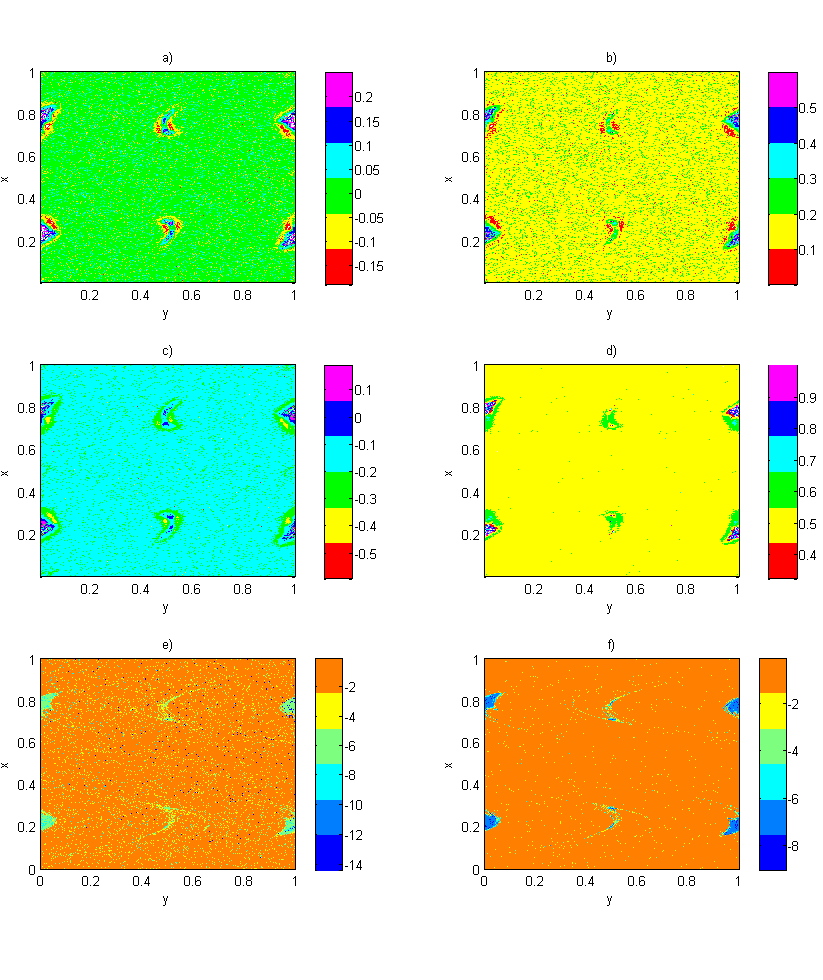}
  \caption{Structure of the Standard Map for $k=4.5$. \textbf{a)}  $\xi'(g_1)$ for $m,n=10^3$, \textbf{b)} $\xi'(g_2)$ for $m, n=10^3$,  \textbf{c)} $\xi'(g_3)$ for $m, n=10^3$, \textbf{d)} $\sigma(g_1)$ for $m, n=10^3$, \textbf{e)} $\log_{10}(R_{t=100})$, \textbf{f)} $\log_{10}(\Delta_{t=100})$}
  \label{65}
\end{figure}

Let us now fix the value of $K=6.5$. We want to show that we can depict the structure of the Standard map with the indicators presented above starting from 500x500 initial conditions uniformly distributed on the bidimensional torus. The number of iterations $k$, $n$, $m$ and $\alpha$ are fixed as before. Results are shown in figure \ref{65} where we compare the four parameters of GEV distribution (top and middle panels) with the Reversibility Error and Divergence of orbits in logarithm scale (lower panel). The number of iterations for the round off indicators is $t=100$. It is evident that the structure of the Standard Map is well highlighted by all the indicators based on GEV distribution. For $g_1$ the empirical values agreed with the theoretical ones $\xi'=0$ and $\sigma=1/2$ expected in the chaotic regions, whereas in the small regular islands we observe significant deviation from the expected values. Similar results hold for $g_2$ and $g_3$ for which the expected theoretical values are $\xi'=1/6$ and $\xi'=-1/6$ respectively. The Round off based indicators highlight the same structure: in the chaotic sea the exponential growth of these quantities lead to a saturation at a value comparable with the size of the torus, whereas in the regular region their values remain order of magnitude lower.\\

\begin{figure}[ht!]
  \includegraphics[width=1.0\textwidth]{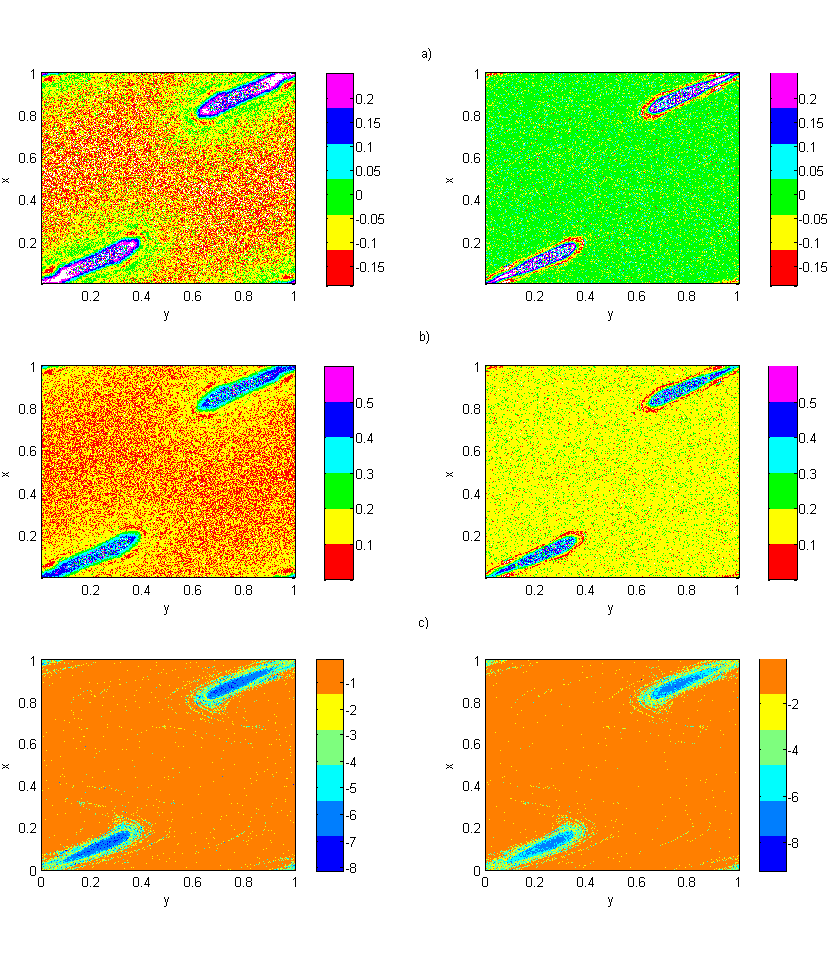}
  \caption{Structure of the Standard Map for $k=4.5$. \textbf{a)} Left-hand side: $\xi'(g_1)$ for $m,n=10^3$, Right-hand side: $\xi'(g_1)$ for $m=10^4, n=10^3$, \textbf{b)} Left-hand side: $\xi'(g_2)$ for $m,n=10^3$, Right-hand side: $\xi'(g_2)$ for $m=10^4, n=10^3$, \textbf{c)} Left-hand side: $\log_{10}(R_{t=100})$, Right-hand side: $\log_{10}(\Delta_{t=100})$ }
  \label{45}
\end{figure}

Eventually, we repeat the same  analysis changing only the parameter $K$ and fixing it to the value $K=4.5$. The structure
of the orbits in the standard map is now quite different, as we encounter larger stability islands and even the motions
that start in the chaotic sea spend time around the regular islands. Due to this reason, we expect to observe a slower convergence
to the parameter of the extreme value distribution with respect to the case presented above.
This depends upon the fact that we may not extract a proper maximum if the motion is affected by a superposition of regular and chaotic motion.
A better convergence must be observed if we increase the observations in each bin to $m=10^4$ as  an extreme value less affected by 
the regular motions can be extracted improving the convergence to the therotical parameters. This is exactly what is shown in figure \ref{45}: in the top panel $\xi'(g_1)$ is compared for the case $m=10^3$ (Left-hand side) and $m=10^4$ (Right-hand side): clearly the convergence
to the expected theoretical values $\xi'=0$ improves when the number of observations in each bin increases. Similar results hold for the
others GEV parameters: as a further example $\xi'(g_2)$ is presented in the middle panel of figure \ref{65}. The structure of the map is again well highlighted by GEV indicators if we compare them with the logarithm of $R_t$, $\Delta_t$ shown in the lower panel for $t=100$.

\section{Conclusions}

Nowadays there exists a large family of dynamical indicators of stability that can be used to obtain informations about the stability
of orbits in dynamical systems. Netherless, the numerical algorithm to compute many of them can be computationally expensive as it usually 
require the solution of variational equations or a Monte Carlo simulation. Moreover each indicator is specialized in detecting
differences either in the chaotic regions or in the regular ones.\\
In this framework we have introduced the parameters of the GEV distribution for some special class of observables to show that they effectively work as other dynamical indicators. As shown in  \citet{freitas}, \citet{faranda2011numerical} and \citet{faranda2011extreme} a GEV like distribution can be observed only if the dynamical system satisfies certain mixing conditions or exhibit an exponential decay of the Hitting Time Statistics.\\
We selected the parameters that do not depend explicitly on the number of observations in each bin $m$, so that our results only depend on the dimension of the attractor $d$. We proved the effectiveness of the indicators testing them on the Standard Map. First we proved that, varying the 
parameters $K$ that regulate the chaoticity of the map, the parameters of the GEV fitted distribution for a small ensemble of initial conditions approach the theoretical values when $K>>1$ that corresponds to chaotic motions. We have also proved that GEV parameters are able to distinguish
regular islands in chaotic sea for Standard Map parameter $K=6.5$ and $K=4.5$. In the latter case we have experienced the presence of superimposition
of regular and chaotic motions that prevented us from obtaining the theoretical expected parameters for the chaotic regions whem $m$ was not big enough.\\
To conclude, the GEV parameters provide basically the same informations of other indicators as the Divergence of 
orbits and the Reversibility error and they are easily accessible from a computational view point as they do not require the solution of the variational equations. With respect to the other indicators they highlight other relevant informations on the dynamics: the parameters are dependent on the dimension of the attractor that can be derived from a fit to the empirical distribution once we are sure that we are dealing
with a chaotic orbit. Furthermore, the extreme value distribution is itself interesting from both practical and theoretical point of view: as far as the observable are concerned it gives detailed informations about the closest return near a certain initial condition.

\section{Aknowledgements}
S.V. was supported by the CNRS-PEPS Project \textit{Mathematical Methods of
Climate Models}, and he thanks the GDRE Grefi-Mefi for having supported
exchanges with Italy. 
DF and VL  acknowledge the financial support of the EU FP7-ERC project NAMASTE Thermodynamics of the Climate System.
DF acknowledges Martin Mestre for the support with the numerical simulations and Peter de Boeck for the useful comments and suggestions that improved the quality of the paper

\bibliographystyle{ws-ijbc}  
\bibliography{indicatori}

\end{document}